\documentclass[a4paper,12pt]{amsart}
\usepackage{amsmath,mathrsfs,amssymb,amsthm,amscd}
\usepackage{a4wide}
\usepackage{hyperref}
\usepackage[all]{xy} 
\usepackage[latin1]{inputenc}

\usepackage{graphicx}
\usepackage{color}
\newcommand{\frakp}{\mathfrak{p}}

\newcommand{\frakq}{\mathfrak{q}}

\newcommand{\CC}{\mathbb{C}}

\newcommand{\PP}{\mathbb{P}}
\newcommand{\QQ}{\mathbb{Q}}
\newcommand{\RR}{\mathbb{R}}
\newcommand{\ZZ}{\mathbb{Z}}

\newcommand{\can}{{\textrm{can}}}

\newcommand{\calC}{{\mathcal{C}}}
\newcommand{\calD}{{\mathcal{D}}}
\newcommand{\calE}{{\mathcal{E}}}
\newcommand{\calF}{{\mathcal{F}}}
\newcommand{\calG}{{\mathcal{G}}}

\newcommand{\calK}{{\mathcal{K}}}
\newcommand{\calL}{{\mathcal{L}}}
\newcommand{\calM}{{\mathcal{M}}}

\newcommand{\calO}{{\mathcal{O}}}

\newcommand{\calS}{{\mathcal{S}}}
\newcommand{\calW}{{\mathcal{W}}}
\newcommand{\calX}{{\mathcal{X}}}
\newcommand{\calY}{{\mathcal{Y}}}
\newcommand{\calV}{{\mathcal{V}}}

\newcommand{\Spec}{\operatorname{Spec}}

\newcommand{\Nm}{\operatorname{Nm}}

\newcommand{\dv}{\operatorname{div}}

\newcommand{\Pic}{\operatorname{Pic}}

\newcommand{\Pica}{\widehat{\operatorname{Pic}}}

\newcommand{\supp}{\operatorname{supp}}
\newcommand{\Ar}{\textrm{Ar}}

\newcommand{\Cl}{\operatorname{Cl}}
\newcommand{\La}[1]{L_{\alpha_{#1}}}
\newcommand{\Laa}{L_{\alpha_{i,j}}}
\newcommand{\dd}{\textrm{d}}

\numberwithin{equation}{section}

\newcommand{\belyi}{{\boldsymbol \beta}}

\newcommand{\GZ}{\mathbb{Z}}
\newcommand{\CZ}{\mathbb{C}}

\newcommand{\fX}{\mathfrak{X}}

\newcommand{\ofin}{\operatorname{fin}}
\newcommand{\finn}[2]{(#1.#2)_{\ofin}}

\newtheorem{satz}{Theorem}[section]
\newtheorem{korollar}[satz]{Corollary}

\newtheorem{lemma}[satz]{Lemma}
\newtheorem{proposition}[satz]{Proposition}

\newenvironment{beweis}{\par\pagebreak[2]\noindent{\bf Proof: }}{
\hfill $\Box$ \bigskip}

\theoremstyle{definition}
\newtheorem{bem}[satz]{Remark}
\newtheorem{defin}[satz]{Definition}

\begin{document}

\title[Arithmetic self-intersection numbers  for Fermat curves of prime exponent]{On the arithmetic self-intersection numbers of 
the dualizing sheaf for Fermat curves of prime exponent}

\author{Christian Curilla, Ulf K\"uhn}
\date{\today}
\email{curilla@math.uni-hamburg.de, kuehn@math.uni-hamburg.de}
\address{Department Mathematik (AZ)\\Universit\"at Hamburg\\
  Bundesstrasse 55\\D-20146 Hamburg}

\begin{abstract}In this article we improve the upper bound for the
  arithmetic self-intersection number of the dualizing sheaf of the
  minimal regular model for the Fermat curves $F_p$ of prime exponent.
\end{abstract}
\maketitle

\tableofcontents

\setcounter{section}{-1}
\section{Introduction}
The main motivation of Arakelov to develop an arithmetic intersection
theory was the idea of proving the Mordell conjecture by mimicking the
proof in the function field case done by Parshin \cite{Par2}. Let $E$ be a number field. A central step
in this program relies on suitable upper bounds for the arithmetic
self-intersection number ${\overline{\omega}_\Ar}^2$,
 where $\overline \omega_{\Ar}$ is the
dualizing sheaf $\omega_{\fX} = \omega_{\fX/\calO_E} \otimes f^*
\omega_{\calO_E/\ZZ}$ equipped with the Arakelov metric (see
\cite{Ar}, p.1177, \cite{MB1}, p.75), 
 of an arithmetic
surface $\fX \to \Spec \calO_E$ that varies in certain complete families (cf. \cite{Par},
\cite{MB2}, or Vojta's appendix in \cite{Lang}). However finding such
bounds turned out to be an intricate problem. 
The best results obtained so far give asymptotics or upper bounds for
${\overline{\omega}_\Ar}^2$ on regular models for certain
discrete families of curves as modular curves (see \cite{A.U},
\cite{M.U}, \cite{JK} and \cite{omega}) and Fermat curves (see \cite{omega}).
Bounds for these curves
have been asked for since the beginning of Arakelov theory 
(see e.g. \cite{Lang}, p. 130 or \cite{MB2}, 8.2).

In this article we improve the upper bound of
${\overline{\omega}_\Ar}^2$ for Fermat curves $F_p$ of prime
exponent. Our calculations rely on a careful analysis of the cusps
behaviour above the prime $p$. This allows us to compute exactly the
``algebraic contributions'' of
a formula for  ${\overline{\omega}_\Ar}^2$ in \cite{omega}.  
We also take into
  account the difference between the minimal regular model
  $\mathfrak{F}_p^{min}$ and the regular model $\mathfrak{F}_p$
  constructed in \cite{Mc}, i.e. the minimal
  desingularisation of the closure in $\PP_{\ZZ [\zeta_p]}^2$ of the
  Fermat curve $x^p+y^p=z^p$ with prime exponent $p$. This leads to
the following result.

\begin{satz}\label{thm:intro}
Let $\pi: \mathfrak{F}_p^{min} \rightarrow \Spec \ZZ[\zeta_p]$ be the minimal regular model of the Fermat curve $F_p:x^p+y^p=z^p$ of prime exponent and genus $g$. Then the arithmetic self-intersection number of its dualizing sheaf equipped with the Arakelov metric satisfies
\begin{equation*} \overline{\omega}_{\mathfrak{F}_p^{min}, \Ar }^2 \le (2 g-2) \left( \log|\Delta_{\QQ(\zeta_p)|\QQ}|^2 + 
   [\QQ(\zeta_p):\QQ]  \left(  \kappa_1 \log p + \kappa_2 \right) + \frac{3p^2-14p+15}{p(p-3)} \log p
   \right),
\end{equation*}
where $\kappa_1, \kappa_2\in \RR$ are  constants independent of $p$.
\end{satz}
It is a well known fact that $\Delta_{\QQ(\zeta_p)|\QQ}=(-1)^{\frac{p-1}{2}}p^{p-2}$ and $[\QQ(\zeta_p):\QQ]=p-1$ and so Theorem \ref{thm:intro} yields 
\begin{equation*} \overline{\omega}_{\mathfrak{F}_p^{min}, \Ar }^2 \le (2 g-2) \left(  
   (p-1)  \left(  \kappa_1 \log p + \kappa_2 \right) + \frac{2p^2-p-5}{p} \log p \right).
\end{equation*}
In comparison to previous results in \cite{omega} our explicit
calculation of the algebraic contributions reduces the maximal
possible growth of $ \overline{\omega}_{\mathfrak{F}_p^{min}, \Ar }^2$
as a function in $p$ by a factor $g(F_p) p^6$.  In the forthcoming
thesis of the first named author the more general case of Fermat
curves with squarefree exponents will be considered.

\section{Intersection theory for arithmetic surfaces}
We start by reminding some notation used in the context of Arakelov Theory. Most of it will be very similar to the notation used in \cite{Soul}. 
\begin{defin}\label{defarithsur} An \emph{arithmetic surface} $\fX$ is a regular integral scheme of dimension
$2$ together with a projective flat morphism $f:\fX \rightarrow
\Spec\calO_{E}$, where $\calO_E$ is the ring of integers of a number
field $E$. Moreover we assume that the \emph{generic fiber}
$X_E=\fX\times_{\Spec\calO_E} \Spec E$ of $f$ is geometrically
irreducible, i.e. $\fX$ is a regular model for $X_E$ over
$\Spec\calO_E$. We denote the complex valued points $\fX(\CZ)$ by
$\fX_\infty$; this is a compact, $1$-dimensional, complex manifold,
which may have several connected components. Actually we have the
decomposition \[ \fX_\infty = \coprod_{\sigma : E \hookrightarrow \CZ
} \fX_\sigma (\CZ) \, ,\] where $\fX_\sigma (\CZ)$ denotes the set of
complex valued points of the curve $\fX_\sigma = \fX \times_{\Spec
  E,\sigma} \Spec\CZ $ coming from the embedding $\sigma
:E\hookrightarrow \CZ$.  For each $s\in \Spec \calO_E$ we define the fibre above $s$ as $\fX_s := \fX \times_{\Spec \calO_E} \Spec k(s)$. We have $\fX_{(0)}=X_E$. Any point $s\neq (0)$ will be called a \emph{closed point} and the corresponding fibre $\fX_s$ a \emph{special fibre}. 
\end{defin}

Let $f:\fX\rightarrow \Spec \calO_E$ be an arithmetic surface in the sense of Definition \ref{defarithsur}. Due to the fact that $\Spec \calO_E$ is Noetherian and that $f$ is of finite type it follows that $\fX$ is Noetherian as well. 

\begin{bem}\label{pic=cl} Since $\fX$ is a regular Noetherian integral scheme, the divisor class group $\Cl (\fX)$ of $\fX$ is isomorphic to the Picard group $\Pic (\fX)$ (see \cite{Liu}, p.257: Corollary 1.19 and p.271: Proposition 2.16). For any divisor $D$ we denote the corresponding invertible sheaf by $\calO_{\fX} (D)$. 
\end{bem}

\begin{defin}
We set $\Cl(\fX)_\QQ =\Cl(\fX)\otimes_\ZZ \QQ$. Obviously $\Cl(\fX)_\QQ$ is a group again. The difference is that we are now allowed to work with divisors with rational coefficients.
\end{defin}

\begin{lemma}\label{hlemma2} Let $f:\fX \rightarrow \Spec \calO_E$ be an arithmetic surface and $s\in \Spec \calO_E$ a closed point. Then \[\fX_s=\frac{1}{m}\dv (h) \] in $\Cl(\fX)_\QQ$, where $\fX_s=f^{\ast} s$, $h\in K(\fX)$ and $m\in \ZZ$.
\end{lemma}

\begin{beweis}
We know that the
  divisor class group $\Cl(\Spec \calO_E)$ is finite and so we can
  find a positive integer $m$ and a rational function $g\in K(\Spec
  \calO_E)$ with the property that $m\cdot s=\dv (g)$. Since $\fX$ is
  regular it follows that $f^{\ast} s=\fX_s$ (see \cite{Liu}, p.351:
  Lemma 3.9) and so $f^{\ast}(m\cdot s)=m\cdot \fX_s=\dv (h)$ in $\Cl(\fX)$ for a
  $h \in K(\fX)$. Now, in $\Cl(\fX)_\QQ$ we may divide this equation by $m$ and the lemma is proven.
\end{beweis}

\begin{defin}
Let $\calD$, $\calE$ be effective divisors without common component, $x\in \fX$ a closed point and $f$, $g$ represent $\calD$, $\calE$ respectively in the local ring $\calO_{\fX,x}$. Then we define the \emph{intersection number} $i_x (\calD,\calE)$ in $x$ as the length of $\calO_{\fX,x}/(f,g)$ as a $\calO_{\fX,x}$-module. The symbol $i_x (\calD,\calE)$ is bilinear and so we may extend the intersection number to all divisors of $\fX$ (just write $\calD$ as $\calD_{+}-\calD_{-}$ with $\calD_{+}$ and $\calD_{-}$ effective and then define $i_x (\calD,\calE):=i_x(\calD_{+},\calE)-i_x(\calD_{-},\calE)$). Now let $s\in\Spec \calO_E$ be a closed point. The \emph{intersection number of $\calD$ and $\calE$ above $s$} is then defined as \[
i_s (\calD,\calE):= \sum_{x\in \fX_s}i_x (\calD,\calE)[k(x)	:k(s)]  \, ,
\] where $x$ runs through the closed points of $\fX_s$ and $k(x)$, $k(s)$ denote the residue class field of $x$, $s$ respectively. If it is clear from the context which intersection number we compute (above which $s$), we simply write $\calD \cdot\calE$.
\end{defin}

\begin{defin}
Let $s\in \Spec\calO_E$ be a closed point and $\calE$ a vertical divisor contained in the special fiber $\fX_s$. According to the moving lemma (see e.g. \cite{Liu}, p.379: Corollary 1.10) there exists a principal divisor $(f)$ so that $\calD:=\calE+(f)$ and $\calE$ have no common component. Since $(f) \cdot \calE=0$ (see. e.g. \cite{Lang}, p.58: Theorem 3.1.) we may define the \emph{self-intersection} of $\calE$ as \[  \calE^2:= \calD \cdot \calE\, .\] 
\end{defin}

\begin{bem}
Another possible way to define $\calE^2$ can be done via cohomological methods (see e.g. \cite{Deli}).
\end{bem}

\section{Canonical divisors on an arithmetic surface}

Let $f:\fX\rightarrow \Spec \calO_E$ be an arithmetic surface in the sense of Definition \ref{defarithsur}.
As $f$ is a local complete intersection (see \cite{Liu}, p.232: Example 3.18.), we can define the canonical sheaf $\omega_{\fX /\Spec \calO_E}$ of $f:\fX \rightarrow \Spec\calO_E$ (see e.g. \cite{Liu}, p.239: Definition 4.7.).

\begin{bem}
Since the scheme $\Spec \calO_E$ is a locally Noetherian scheme and $f$ is a flat projective local complete intersection of relative dimension $1$, the canonical sheaf is isomorphic to the 1-dualizing sheaf (see \cite{Liu}, p.247: Theorem 4.32.).
\end{bem}

\begin{defin} We call any divisor $\calK$ of $\fX$ with $\calO_{\fX} (\calK) \cong \omega_{\fX /\Spec \calO_E}$ a \emph{canonical divisor}. This divisor exists because of Remark \ref{pic=cl}.
\end{defin}

\begin{bem} Let $s\in \Spec \calO_E$ be a closed or the generic point. For each fibre $\fX_s \rightarrow \Spec k(s)$ we get a canonical sheaf $\omega_{\fX_s / \Spec k(s)}$. We have the relation $\omega_{\fX_s / \Spec k(s)} \cong \omega_{\fX/ \Spec \calO_E} |_{\fX_s} $ (see \cite{Liu}, p.239: Theorem 4.9). If $s$ is the generic point we can define a canonical divisor $K$ of $X:=\fX \times_{\Spec \calO_E} \Spec E$ in the same way we did with the arithmetic surface. Similar to the relation between the canonical sheaves we get $\calK |_{X} \cong K$.
\end{bem}

Now let $\calE$ be a vertical divisor contained in a special fiber $\fX_s$ and $\calK$ a canonical divisor on $\fX$. Since any other canonical divisor is rationally equivalent to $\calK$ the intersection number $\calK \cdot \calE$ depends uniquely on $\omega_{\fX /\calO_{\Spec E}}$ and not on the choice of a representative $\calK$. We have the following important theorem:

\begin{satz}[Adjunction formula]
Let $f:\fX \rightarrow
\Spec\calO_{E}$ be an arithmetic surface, $s\in \Spec \calO_E$ a closed point and $\calE$ a vertical divisor contained in the special fiber $\fX_s$. Then we have \begin{equation}\label{adform} 2p_a(\calE)-2=\calE^2+\calK\cdot \calE \, ,\end{equation}
where $p_a(\calE)$ is the arithmetic genus of $\calE$.
\end{satz}
\begin{beweis}
See \cite{Licht} Theorem 3.2.
\end{beweis}

Later on it will be important to construct the canonical divisor explicitly. The following proposition will help us with that.

\begin{proposition}\label{hprop1}
Let $\calC\in \Cl_{\mathbb{Q}}(\fX)$ be a divisor on $\fX$ which satisfies the adjunction formula \eqref{adform} and whose restriction to the generic fibre $X$ is a canonical divisor of $X$. Then $\calC$ is a canonical divisor on $\fX$.
\end{proposition}

\begin{beweis}
  Let $\calK$ be a canonical divisor on $\fX$ (we already know that it
  exists). We want to show that $\calK \sim \calC$ and so that
  $\calC$ is a canonical divisor as well. We denote the horizontal
  part of the divisors by $\calK_h$ and $\calC_h$. Since the
  restriction to the generic fibre of both divisors is a canonical
  divisor of $X$ we have $\calK |_{X}=\calK_h |_{X}\sim \calC_h
  |_{X}=\calC |_{X}$ and so there exists a rational element $g\in
  K(X)$, which yields $\calK |_{X}-\dv (g)=\calC |_{X}$. Because we
  have $K(X)\cong K(\fX)$, we can interpret $g$ as an element of
  $K(\fX)$ and so obtain a
  principal divisor whose restriction to $X$ is $\dv (g)$. We denote
  this principal divisor by $\dv (g)$ as well. If we now set
  $\calC^{'}:=\calC+\dv (g)$ we get a divisor with the properties that
  $\calC^{'} \sim \calC$ and $\calC^{'}_h=\calK_h$. Since we are
  just interested in $\calC$ up to rational equivalence we may assume
  from now on that the horizontal part of $\calC$ is the same as the
  one of $\calK$.\\ Let $s\in \Spec \calO_E$ be a closed point and
  $\fX_s$ the fibre above it. We denote by $\calK_s$ and $\calC_s$ the
  vertical divisor of $\calK$ and $\calC$ which have support in
  $\fX_s$. Since $\calK$ and $\calC$ fulfill the adjunction formula and have the same horizontal part 
  we have \[
  0=(\calK_s-\calC_s)\cdot(\calK-\calC)=(\calK_s-\calC_s)\cdot(\calK_s-\calC_s)
  \, .\]and so $\calK_s-\calC_s = q\fX_s$, where $q$ is a rational
  number (see \cite{Lang}, p.61: Proposition 3.5.). Now, according to Lemma \ref{hlemma2}, we find $m\in\ZZ$ and $h\in K(\fX)$ so that $\calK_s-\calC_s =q\fX_s =\frac{q}{m} \dv
  (h)$ and so we have $\calK_s \sim \calC_s$ in
  $\Cl(\fX)_{\mathbb{Q}}$. If we set $\calC^{'}:=\calC+\frac{q}{m}
  \dv (h)$ we have just changed the components of $\calC$ with
  support in $\fX_s$. Again, we have $\calC^{'}\sim \calC$ and now
  $\calK_h+\calK_s = \calC^{'}_h +\calC^{'}_s$. Continuing
  successively with the other closed points of $\Spec \calO_E$ we
  arrive at a divisor $\calC^{''}$ with $\calC^{''}=\calK$ and
  $\calC^{''}\sim \calC$ as we claimed at the beginning.
\end{beweis}

\begin{bem}\label{bem:rational} The Proposition \ref{hprop1} uses the fact that in $\Cl(\fX)_\QQ$ the special fibres are divisors coming from functions (see Lemma \ref{hlemma2}). In other words, the canonical divisor is only defined up to rational multiples of the special fibres (in $\Cl(\fX)_\QQ$).

\end{bem}

\section{Arithmetic intersection numbers for hermitian line bundles}
\begin{defin}
A \emph{hermitian line bundle} $\overline{\calL}=(\calL,h)$ is a line bundle $\calL$ on $\fX$ together with a smooth, hermitian metric $h$ on the induced holomorphic line bundle $\calL_\infty =\calL \otimes_\GZ \CZ$ on $\fX_\infty$. We denote the norm associated with $h$ by $||\cdot ||$. Two hermitian line bundles $\overline{\calL}$, $\overline{\calM}$ on $\fX$ are \emph{isomorphic}, if \[ \overline{\calL}\otimes \overline{\calM}^{^{-1}} \cong (\calO_{\fX}, |\cdot|) \, ,\] where $|\cdot|$ denotes the usual absolute value. The \emph{arithmetic Picard group} $\Pica (\fX)$ is the group of isomorphy classes of hermitian line bundles $\overline{\calL}$ on $\fX$, the group structure being given by the tensor product.
\end{defin} 

\begin{defin}
Let $\overline{\calL}$, $\overline{\calM}$ be two hermitian line bundles on $\fX$ and $l,m$ non-trivial, global sections, whose induced divisors $\dv (l)$ and $\dv (m)$ on $\fX$ have no horizontal component in common. Then we define the \emph{intersection number at the finite places} $\finn{l}{m}$ of $l$ and $m$ by the formula \begin{eqnarray*}
\finn{l}{m} &:=&\sum_{x\in\fX} \log\sharp \left(\calO_{\fX,x} / (l_x,m_x) \right) = \sum_{x\in \fX} i_x (\dv (l),\dv (m)) \log \left|k(x)\right| \\ &=&\sum_{s\in \Spec \calO_E} \left(\sum_{x\in \fX_s}i_x (\dv (l),\dv (m))[k(x)	:k(s)] \right)\log \left|k(s)\right|\, ,
\end{eqnarray*}where $l_x$ and $m_x$ are local equations of $l$ and $m$ at the point $x\in \fX$; the sum runs through the closed points $x$ of $\fX$. \\ The sections $l$ and $m$ induce global sections on $\calL_\infty$ and $\calM_\infty$, which we denote by abuse of notation again by $l$ and $m$. We assume that the associated divisors $\dv (l)$ and $\dv (m)$ on $\fX_\infty$ have no points in common. Writing $\dv (l) =\sum_{\alpha} p_\alpha P_\alpha$ with $p_\alpha \in \GZ$ and $P_\alpha \in \fX_\infty$, we set \[ (\log ||m||)[\dv (l)]:= \sum_{\alpha} p_\alpha \log ||m(P_\alpha)|| \, .\] The \emph{intersection number at the infinite places} $(l.m)_\infty$ of $l$ and $m$ is now given by the formula \[
(l.m)_\infty := -(\log ||m||)[\dv (l)]-\int_{\fX_\infty} \log ||l||\cdot c_1 (\overline{\calM}) \, ,
\] where the first Chern form $c_1 (\overline{\calM})\in H^{1,1}(\fX_\infty, \RR)$ of $\overline{\calM}$ is given, away from the divisor $\dv (m)$ on $\fX_\infty$, by \[ c_1 (\overline{\calM})=\operatorname{dd^c}(-\log ||m(\cdot)||^2) \, .\]
We define the \emph{arithmetic intersection number} $\overline{\calL}.\overline{\calM}$ of $\overline{\calL}$ and $\overline{\calM}$ by \begin{equation}\label{arak}
\overline{\calL}.\overline{\calM}:=\finn{l}{m}+(l.m)_\infty \, .
\end{equation}The \emph{arithmetic self-intersection number} of $\overline{\calL}$ is given by $\overline{\calL}.\overline{\calL}$. 

\end{defin}
\begin{satz}[Arakelov, Deligne et al.]\label{aip} Formula \eqref{arak} induces a bilinear, symmetric pairing \[ \Pica(\fX)\times\Pica(\fX)\rightarrow \RR \, .\]
\end{satz}
\begin{beweis}
See for example \cite{Soul}.
\end{beweis}

\begin{bem}
Theorem \ref{aip} is a generalisation, essentially due to Deligne, of the arithmetic intersection pairing, invented by Arakelov, where only hermitian line bundle, whose Chern forms are multiples of a fixed volume form, are considered.
\end{bem}

If the genus of $\fX$ is greater than one, then for each $\sigma$ we
have on $\fX_\sigma(\CC)$ the \emph{canonical volume form}
\begin{align*}
\nu_\can^\sigma (z) = \frac{i}{2g} \sum_j |f_j^\sigma|^2 \dd z \land \dd \overline z,
\end{align*}
where $f_1^\sigma(z) \dd z$, ... $f_g^\sigma(z) \dd z$ is an
orthonormal basis of $H^0(\fX_\sigma(\CC),\Omega^1)$ equipped with the
natural scalar product. We write $\nu_\can$ for the induced volume
form on $\fX_\infty$ and for ease of notation we set 
$$
\overline\calO(D) = \overline\calO(D)_{\nu_\can}.
$$
Here the norm of the section $1_D$ of $\calO (D)$ is given by $\| 1_D \|=g(D, \cdot)$ where $g$ is the canonical green function (see e.g. \cite{Lang}).

 Due to Arakelov is the observation that there is a unique metric
  $\|\cdot \|_\Ar$ on $\omega_{\fX}$ such that for all sections $P$
  of $\fX$ it holds the adjunction formula
\begin{align} \label{eq:adjunction}
\overline{\omega}_\Ar. \overline{\calO}(P)  + 
\overline{\calO}(P) ^2 = \log| \Delta_{E|\QQ}|, 
\end{align} where $\overline{\omega}_\Ar =(\omega_{\fX},\|\cdot \|_\Ar)$.
Moreover $\overline{\omega}_\Ar$ is a
$\nu_\can$-admissible line bundle (see \cite{Lang}).  
\begin{bem}
In Remark \ref{bem:rational} we saw that the canonical divisor is only defined up to rational multiples of the special fibres. Because of formula \eqref{eq:adjunction} this indeterminacy will be deleted by the norm of the section.
\end{bem}


Let $\calY \to \Spec \calO_E$ be an arithmetic surface and 
write $Y$ for its generic fiber. We fix
$\infty , P_1 , ...,P_r \in Y(E)$ such that $Y \setminus
\{\infty,P_1, ...,P_r\}$ is hyperbolic. Then we
consider any arithmetic surface $\calX\to \Spec \calO_E$ 
equipped with a morphism of arithmetic surfaces 
$\belyi: \calX \to \calY$ such that the induced  morphism
$\belyi: X \to Y$  of algebraic curves defined over  $E$ is  unramified above
$Y(E) \setminus\{\infty ,P_1,...,P_r \}$. 
Let $g\ge 2 $ be the genus of $X$ and $d = \deg(\belyi)$. We write $\belyi^* \infty = \sum b_j S_j$ and 
the points $S_j$ will be called  cusps. Set $b_{\max} = \max_j \{b_j\}$. 
Divisors on $X$ with support in the cusps of degree zero are called cuspidal. Finally, a prime $\frakp$ is said to be bad if the fiber of $\calX$ above $\frakp$ is reducible\footnote{note that a prime of bad reduction need not 
be a bad prime}.

\begin{satz} \label{thm:keyformula} Let $\belyi: \calX \to \calY$ be 
a morphism of arithmetic surfaces as above.
Assume that all cusps are
  $E$-rational points and  that all cuspidal divisors  are torsion, then
  the arithmetic self-intersection number of the dualizing sheaf on  $\calX$ 
  satisfies the inequality
\begin{align}\label{eq:main}
  \overline{\omega}_\Ar^2 &\le (2 g-2) \left( \log|\Delta_{E|\QQ}|^2 + 
   [E:\QQ]  \left(  \kappa_1 \log b_{\max} + \kappa_2 \right) + \sum_{\frakp \,\,{\rm bad}} a_\frakp
  \log\Nm(\frakp)
   \right),
\end{align}
where $\kappa_1, \kappa_2\in \RR$ are  constants
that dependent only on $Y$ and the points $\infty,P_1,...,P_r$.
The coefficients $a_\frakp \in \QQ$
are determined by certain local intersection numbers (see formula \eqref{eq:def-ap} below).
\end{satz}

\begin{beweis}
See \cite{omega} Theorem I.
The method of proof uses classical Arakelov theory, as well as generalized 
arithmetic intersection theory (see \cite{gain}), which allows to use results of Jorgenson and Kramer \cite{JK2}.
\end{beweis}

To keep the notation simple, we write $S_j$ also for the Zariski
closure in $\calX$ of a cusp $S_j$.  Let $\calK$ be a canonical
divisor of $\calX$, then for each cusp $S_j$ we can find a divisor
$\calF_j$ such that
\begin{align*}
 \left( S_j + \calF_j -\frac{1}{2g-2} \calK\right )\cdot\mathcal{C}^{(\mathfrak{p})}_{l}=0
\end{align*}
for all irreducible components $\mathcal{C}^{(\mathfrak{p})}_{l}$ of
the fiber $f^{-1}(\mathfrak{p})$ above $\mathfrak{p} \in \Spec
\mathcal{O}_K$. 
Similarly we find  
 for each cusp $S_j$   a
divisor $\calG_j$ such that also for all $ \mathcal{C}^{(\mathfrak{p})}_{l}$ as before
\begin{align*}
\left ( S_j + \calG_j -\frac{1}{d} \belyi^*\infty \right) \cdot \mathcal{C}^{(\mathfrak{p})}_{l}=0.
\end{align*}

Then the rational numbers $a_\frakp$ in the theorem are determined by 
the following arithmetic intersection numbers
of trivially metrised hermitian line bundles 
\begin{align}\label{eq:def-ap}
 \sum_{\frakp \,\,{\rm bad}} a_\frakp
  \log\Nm(\frakp)
&= - \frac{2 g 
}{d} \sum_j b_j\,\mathcal{O}( \calG_j)^2 +
\frac{2g-2}{d} \sum_j b_j\, \calO(\calF_j)^2.   
\end{align}

\section{Fermat curves and their natural Belyi uniformization}\label{sec:belyimorph}

For the rest of this article we will consider the Fermat curve
\begin{align*}
F_p: X^p+Y^p=Z^p,
\end{align*}
where $p>3$ is prime number, together with the natural
 morphism 
\begin{align}\label{derMor}
\beta:F_p\rightarrow \mathbb{P}^1
\end{align}
 given by
$(x:y:z)\mapsto (x^p:y^p)$.  Since the morphism $\beta$ is defined over $\QQ$, it is defined
 over any number field.
It is a Galois covering of degree $p^2$ and, since there are
only the three branch points $0,1,\infty$, it is a Belyi morphism. 
All the ramification orders equal $p$.
In \cite{MurtyRama} Murty and Ramakrishnan give the associated Belyi uniformisation 
$F_p(\mathbb{C})\setminus \belyi^{-1}\{0,1,\infty\} \cong \Gamma_P 
\setminus \mathbb{H}$. The subgroup $\Gamma_P$
of $\Gamma(2)$ is given by $\Gamma_p = \ker \psi$ where $\psi:\Gamma(2) \to \GZ / p\GZ \times \GZ / p\GZ$ maps the generators of $\Gamma(2)$ to the elements
$(1,0)$ and $(0,1)$. 

A
ramified point, i.e. an element $S\in F_p$ that maps to one of the
branch points, will be called a \emph{cusp}. Divisors with support in the cusps having degree zero are called \emph{cuspidal divisor}.

\begin{proposition}\label{prop:Rohrlich}
Let $F_p$ a Fermat curve and $\beta:F_p\rightarrow \mathbb{P}^1$ the morphism in \eqref{derMor}. \begin{enumerate}
\item The group of cuspidal divisors is a torsion subgroup of $\Cl(F_p)$.
\item Let $S\in F_p(\QQ (\zeta_p))$ be a cusp, then $(2g-2)S$ is a canonical divisor. 
\end{enumerate}
\end{proposition}

\begin{beweis}
The first statement follows from \cite{Roh}, p. 101: Theorem 1. So only the second statement is left. By the Hurwitz formula there exists a canonical divisor with support in the cusps. Then by (i) the claim follows.
\end{beweis}
\section{A regular model and the minimal model for \texorpdfstring{$F_p$}{Fp}}
\label{regandminmodel}
\begin{figure}[h]\centering
\begin{picture}(0,0)%
\includegraphics{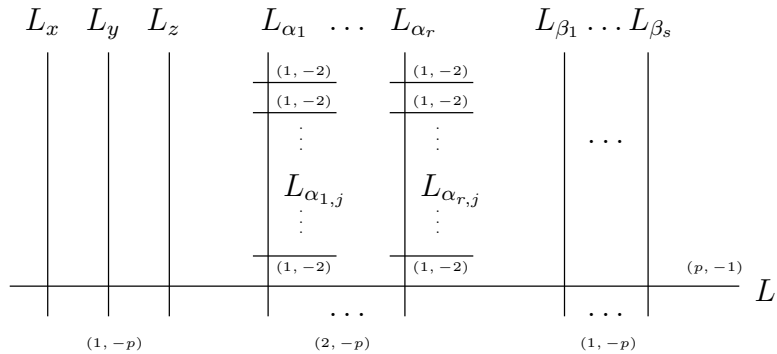}%
\end{picture}%
\setlength{\unitlength}{4144sp}%
\begingroup\makeatletter\ifx\SetFigFontNFSS\undefined%
\gdef\SetFigFontNFSS#1#2#3#4#5{%
  \reset@font\fontsize{#1}{#2pt}%
  \fontfamily{#3}\fontseries{#4}\fontshape{#5}%
  \selectfont}%
\fi\endgroup%
\begin{picture}(4437,2115)(5164,-5638)
\put(5986,-3706){\makebox(0,0)[lb]{\smash{{\SetFigFontNFSS{12}{14.4}{\familydefault}{\mddefault}{\updefault}{\color[rgb]{0,0,0}$L_z$}%
}}}}
\put(5266,-3706){\makebox(0,0)[lb]{\smash{{\SetFigFontNFSS{12}{14.4}{\familydefault}{\mddefault}{\updefault}{\color[rgb]{0,0,0}$L_x$}%
}}}}
\put(5626,-3706){\makebox(0,0)[lb]{\smash{{\SetFigFontNFSS{12}{14.4}{\familydefault}{\mddefault}{\updefault}{\color[rgb]{0,0,0}$L_y$}%
}}}}
\put(7066,-5416){\makebox(0,0)[lb]{\smash{{\SetFigFontNFSS{12}{14.4}{\familydefault}{\mddefault}{\updefault}{\color[rgb]{0,0,0}$\hdots$}%
}}}}
\put(9586,-5326){\makebox(0,0)[lb]{\smash{{\SetFigFontNFSS{12}{14.4}{\familydefault}{\mddefault}{\updefault}{\color[rgb]{0,0,0}$L$}%
}}}}
\put(6796,-4696){\makebox(0,0)[lb]{\smash{{\SetFigFontNFSS{12}{14.4}{\familydefault}{\mddefault}{\updefault}{\color[rgb]{0,0,0}$L_{\alpha_{1,j}}$ }%
}}}}
\put(7606,-4696){\makebox(0,0)[lb]{\smash{{\SetFigFontNFSS{12}{14.4}{\familydefault}{\mddefault}{\updefault}{\color[rgb]{0,0,0}$L_{\alpha_{r,j}}$ }%
}}}}
\put(7696,-4426){\makebox(0,0)[lb]{\smash{{\SetFigFontNFSS{5}{6.0}{\familydefault}{\mddefault}{\updefault}{\color[rgb]{0,0,0}$\vdots$}%
}}}}
\put(6886,-4426){\makebox(0,0)[lb]{\smash{{\SetFigFontNFSS{5}{6.0}{\familydefault}{\mddefault}{\updefault}{\color[rgb]{0,0,0}$\vdots$}%
}}}}
\put(6886,-4921){\makebox(0,0)[lb]{\smash{{\SetFigFontNFSS{5}{6.0}{\familydefault}{\mddefault}{\updefault}{\color[rgb]{0,0,0}$\vdots$}%
}}}}
\put(7696,-4921){\makebox(0,0)[lb]{\smash{{\SetFigFontNFSS{5}{6.0}{\familydefault}{\mddefault}{\updefault}{\color[rgb]{0,0,0}$\vdots$}%
}}}}
\put(6751,-5146){\makebox(0,0)[lb]{\smash{{\SetFigFontNFSS{5}{6.0}{\familydefault}{\mddefault}{\updefault}{\color[rgb]{0,0,0}$(1,-2)$}%
}}}}
\put(7561,-5146){\makebox(0,0)[lb]{\smash{{\SetFigFontNFSS{5}{6.0}{\familydefault}{\mddefault}{\updefault}{\color[rgb]{0,0,0}$(1,-2)$}%
}}}}
\put(7561,-4156){\makebox(0,0)[lb]{\smash{{\SetFigFontNFSS{5}{6.0}{\familydefault}{\mddefault}{\updefault}{\color[rgb]{0,0,0}$(1,-2)$}%
}}}}
\put(7561,-3976){\makebox(0,0)[lb]{\smash{{\SetFigFontNFSS{5}{6.0}{\familydefault}{\mddefault}{\updefault}{\color[rgb]{0,0,0}$(1,-2)$}%
}}}}
\put(6751,-4156){\makebox(0,0)[lb]{\smash{{\SetFigFontNFSS{5}{6.0}{\familydefault}{\mddefault}{\updefault}{\color[rgb]{0,0,0}$(1,-2)$}%
}}}}
\put(6751,-3976){\makebox(0,0)[lb]{\smash{{\SetFigFontNFSS{5}{6.0}{\familydefault}{\mddefault}{\updefault}{\color[rgb]{0,0,0}$(1,-2)$}%
}}}}
\put(5626,-5596){\makebox(0,0)[lb]{\smash{{\SetFigFontNFSS{5}{6.0}{\familydefault}{\mddefault}{\updefault}{\color[rgb]{0,0,0}$(1,-p)$}%
}}}}
\put(8281,-3706){\makebox(0,0)[lb]{\smash{{\SetFigFontNFSS{12}{14.4}{\familydefault}{\mddefault}{\updefault}{\color[rgb]{0,0,0}$L_{\beta_1}\ldots L_{\beta_s}$}%
}}}}
\put(9181,-5146){\makebox(0,0)[lb]{\smash{{\SetFigFontNFSS{5}{6.0}{\familydefault}{\mddefault}{\updefault}{\color[rgb]{0,0,0}$(p,-1)$}%
}}}}
\put(8596,-4381){\makebox(0,0)[lb]{\smash{{\SetFigFontNFSS{12}{14.4}{\familydefault}{\mddefault}{\updefault}{\color[rgb]{0,0,0}$\hdots$}%
}}}}
\put(8596,-5416){\makebox(0,0)[lb]{\smash{{\SetFigFontNFSS{12}{14.4}{\familydefault}{\mddefault}{\updefault}{\color[rgb]{0,0,0}$\hdots$}%
}}}}
\put(8551,-5596){\makebox(0,0)[lb]{\smash{{\SetFigFontNFSS{5}{6.0}{\familydefault}{\mddefault}{\updefault}{\color[rgb]{0,0,0}$(1,-p)$}%
}}}}
\put(6976,-5596){\makebox(0,0)[lb]{\smash{{\SetFigFontNFSS{5}{6.0}{\familydefault}{\mddefault}{\updefault}{\color[rgb]{0,0,0}$(2,-p)$}%
}}}}
\put(6661,-3706){\makebox(0,0)[lb]{\smash{{\SetFigFontNFSS{12}{14.4}{\familydefault}{\mddefault}{\updefault}{\color[rgb]{0,0,0}$L_{\alpha_1}$ }%
}}}}
\put(7426,-3706){\makebox(0,0)[lb]{\smash{{\SetFigFontNFSS{12}{14.4}{\familydefault}{\mddefault}{\updefault}{\color[rgb]{0,0,0}$L_{\alpha_r}$ }%
}}}}
\put(7066,-3706){\makebox(0,0)[lb]{\smash{{\SetFigFontNFSS{12}{14.4}{\familydefault}{\mddefault}{\updefault}{\color[rgb]{0,0,0}$\hdots$}%
}}}}
\end{picture}%

\caption{The configuration of the geometric special fibre $\mathfrak{F}_{p} \times_{\Spec R}\Spec{\overline{k(\pi)}}$.}\label{d1}
\end{figure}

In this section we are going to sketch the construction done by
McCallum \cite{Mc} of a regular model and the minimal model of the
curve $ {F_p}: x^p+y^p=z^p$ over $S=\Spec R$, where
$R=\ZZ_p [\zeta_p]$ denotes the ring of integers of the field $\mathbb{Q}_p (\zeta_p)$ and $\zeta_p$ a primitive $p$-th root of
unity. In order to simplify our computations we may consider the curve
\begin{equation}\label{1} C_p : x^p+y^p=1 \end{equation} in
$\mathbb{A}^2_S$ because the model, we are starting with, is just the
normalization of the projective completion of $C_p$. It has just one
prime ideal of bad reduction, namely $(\pi):=(1-\zeta_p)$ which is the
only prime lying over $(p)$; in fact since $p$ is totally ramified in
$\mathbb{Q}_p (\zeta_p)$ we have $p=u\pi^{p-1}$ with an element $u\in
\mathbb{Z}_p [\zeta_p]^{\ast}$. Reduction modulo $p$ gives us a $p$-tuple
line which is non-regular. Moving this line to the $x$-axis, or in other words setting \begin{equation}\label{transf}X=x~ \mbox{   and   }  ~Y=y+x-1 \,, \end{equation} equation
\eqref{1} becomes \[ -u\pi^{p-1}\phi(X,-Y-1)+u\pi^{p-1}\phi(Y)+Y^p=0 \,
,\] where \[ \phi(X,Y):=\frac{(X+Y)^p-X^p-Y^p}{p}
\] and $\phi(X):=\phi(X,1)$. Now, by blowing up the line $\pi=Y=0$, one obtains a model which is covered by the two affine open sets: we introduce new variables $a$ and $b$. Setting $b=\frac{\pi}{Y}$, we have $U_1=\Spec \left(R[X,Y,b]/(bY-\pi,F_1(X,Y) )\right)$ where \[ F_1(X,Y)= -u b^{p-1}\phi(X,-Y-1)+u b^{p-1}\phi(Y)+Y\, ; \] setting $a=\frac{Y}{\pi}$ the second affine open set is $ U_2=\Spec \left(R[X,Y,a]/(a\pi-Y,F_2(X,Y)) \right) $ where \[ F_2(X,Y)=-u\phi(X,-Y-1)+u\phi(Y)+\pi a^p \,.\] The geometric special fibre $U_1\times_S \Spec \overline{k(\pi)}\cup U_2\times_S \Spec \overline{k(\pi)}$ of this model consists of a
component $L$ (which is located just in $U_1$ and associated to the ideal $(\overline{Y},\overline{b})$ in $R[X,Y,b]/(bY-\pi,F_1(X,Y)$) and components $L_x,  L_y
$, $L_{\alpha_1},\ldots,L_{\alpha_r}, L_{\beta_1},\ldots,L_{\beta_s}$
which intersect $L$ and correspond to the different roots of the
polynomial \[ \phi(X,-1)=-X(X-1)\prod_{\genfrac{}{}{0pt}{}{\alpha \neq 0,1}{  \alpha
  \in k(\pi)}}(X-\alpha)^2 \prod_{\beta \notin k(\pi)}(X-\beta) \, .\]
The $L_{\alpha_i}$ appear with multiplicity $2$ whereas all other
components with multiplicity $1$. There is also a line $L_z$ crossing
the point at infinity on $L$, which we cannot see in this affine
model. There are just singularities left on the double lines
$L_{\alpha_i}$. Blowing up these singularities we achieve new
components $L_{\alpha_{i,j}}$ crossing $L_{\alpha_{i}}$. All
components have genus 0. For later applications we define the index set \begin{equation}\label{indexset}I:=\left\{x,y,z,\beta_i,\alpha_j, \alpha_{j,k},\ldots \right\} .\end{equation} Let us denote the model we achived by
$\mathfrak{F}_p$. The scheme $\mathfrak{F}_p$ is a regular model and
its geometric special fibre $\mathfrak{F}_{p} \times_{\Spec R}\Spec{\overline{k(\pi)}}$
corresponding to $(\pi)$ has the configuration as in figure \ref{d1}
where all components of the fibre have genus $0$ and the pair
$(n,m)$ indicates the multiplicity $n$ and the self-intersection $m$
of the component (\cite{Mc}, Theorem 3.).

\begin{bem} If we now blow down the curve $L$ (which is the only one with self-intersection $-1$), we get the minimal regular model $\mathfrak{F}^{min}_p $ (see \cite{Chin}, p.315: Theorem 3.1). 
\end{bem}

\begin{bem}
A regular model over $\ZZ[\zeta_p]$ can be obtained by glueing the model $\mathfrak{F}_p$ over $S$ and the smooth model of $F_p$ over $\Spec \ZZ[\zeta_p] \setminus \left\{\frakp\right\}$. We will denote this model as well by $\mathfrak{F}_p$.
\end{bem}

Since we were just performing a sequence of blow-ups, the morphism
$\beta:F_{p}\rightarrow \mathbb{P}^1$ extends to a morphism of
arithmetic surfaces
\begin{equation*}
  \beta:\mathfrak{F}_p: \rightarrow \mathbb{P}_{\ZZ[\zeta_p]}^{1}. 
\end{equation*} 
In particular together with Proposition \ref{prop:Rohrlich} we see
that $\beta$ fulfills the assumptions of Theorem
\ref{thm:keyformula}. The rest of this paper is devoted to calculate
the quantities $a_\frakp$ in this theorem.


\section{Extensions of cusps and canonical divisors on \texorpdfstring{$\mathfrak{F}_p$}{Fp}}
\begin{defin} We denote by $S_x$ a cusp of the form $(0:\zeta_p^{i}:1)$; this abuse of notation will be
justified by the Lemma \ref{justi} below, which shows that the
properties of $S_x$, relevant for our considerations, do not depent on
the exponent $i$. Similar we denote by $S_y$ (resp. $S_z$) a cusp of
the form $(\zeta_p^{i}:0:1)$ (resp. $(\zeta_p^{i}:-1:0)$). If we
take the Zariski-closure of a cusp $S_x$ in $\mathfrak{F}_p$, we get a
horizontal divisor, which we denote by $\calS_x$. Again, similar for
$y$ and $z$. \end{defin}

For any two divisors $D$ and $E$ of $\mathfrak{F}_p$ we
say that $D$ \emph{intersects} $E$, if $\supp D \cap \supp E \neq 0$.
 
\begin{figure}[htb]\centering
\begin{picture}(0,0)%
\includegraphics{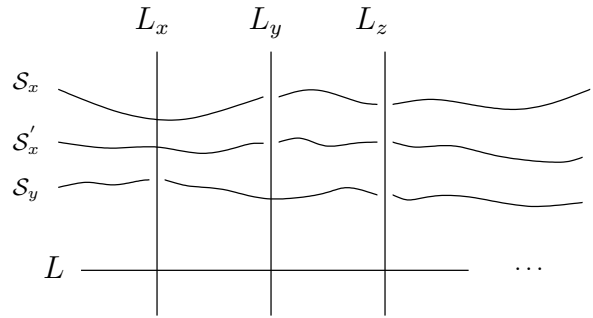}%
\end{picture}%
\setlength{\unitlength}{4144sp}%
\begingroup\makeatletter\ifx\SetFigFontNFSS\undefined%
\gdef\SetFigFontNFSS#1#2#3#4#5{%
  \reset@font\fontsize{#1}{#2pt}%
  \fontfamily{#3}\fontseries{#4}\fontshape{#5}%
  \selectfont}%
\fi\endgroup%
\begin{picture}(3447,1905)(4531,-5428)
\put(5266,-3706){\makebox(0,0)[lb]{\smash{{\SetFigFontNFSS{12}{14.4}{\familydefault}{\mddefault}{\updefault}{\color[rgb]{0,0,0}$L_x$}%
}}}}
\put(6571,-3706){\makebox(0,0)[lb]{\smash{{\SetFigFontNFSS{12}{14.4}{\familydefault}{\mddefault}{\updefault}{\color[rgb]{0,0,0}$L_z$}%
}}}}
\put(5941,-3706){\makebox(0,0)[lb]{\smash{{\SetFigFontNFSS{12}{14.4}{\familydefault}{\mddefault}{\updefault}{\color[rgb]{0,0,0}$L_y$}%
}}}}
\put(7516,-5146){\makebox(0,0)[lb]{\smash{{\SetFigFontNFSS{9}{10.8}{\familydefault}{\mddefault}{\updefault}{\color[rgb]{0,0,0}$\hdots$}%
}}}}
\put(4726,-5191){\makebox(0,0)[lb]{\smash{{\SetFigFontNFSS{12}{14.4}{\familydefault}{\mddefault}{\updefault}{\color[rgb]{0,0,0}$L$}%
}}}}
\put(4546,-4066){\makebox(0,0)[lb]{\smash{{\SetFigFontNFSS{9}{10.8}{\familydefault}{\mddefault}{\updefault}{\color[rgb]{0,0,0}$\calS_x$}%
}}}}
\put(4546,-4426){\makebox(0,0)[lb]{\smash{{\SetFigFontNFSS{9}{10.8}{\familydefault}{\mddefault}{\updefault}{\color[rgb]{0,0,0}$\calS_x^{'}$}%
}}}}
\put(4546,-4696){\makebox(0,0)[lb]{\smash{{\SetFigFontNFSS{9}{10.8}{\familydefault}{\mddefault}{\updefault}{\color[rgb]{0,0,0}$\calS_y$}%
}}}}
\end{picture}%
\caption{The divisors $\calS_x,\calS_x^{'}$ and $\calS_y$, where $\calS_x^{'}$ is coming from another cusp of the form $(0:\zeta_p^{j}:1)$.}\label{d2}
\end{figure}

\begin{proposition}\label{justi}
  Let $\calS$ and $\calS^{'}$ be horizontal divisors of
  $\mathfrak{F}_p$ coming from different cusps $S$ and $S^{'}$ on
  $F_p$. Then the following properties are true:
  \begin{enumerate}
\item $\calS$ does not intersect $\calS^{'}$. 
\item If $\calS=\calS_x$ (resp. $\calS_y,\calS_z$), then $\calS$ only
    intersects the component $L_x$ (resp. $L_y,L_z$) in the special
    fiber $\mathfrak{F}_{p} \times_{\Spec \calO_e}\Spec{k(\pi)}$ (see figure
    \ref{d2}).
\end{enumerate}
\end{proposition}

\begin{beweis}
For the proof, we need to work with the explicit description of the regular model $\mathfrak{F}_p$. So if we talk about a cusp in the following, we will mean a point of the form $(0:\zeta_p^{i}-1:1)$ ($( \zeta_p^{i}:\zeta_p^{i}-1:1)$ resp.) which is just $S_x$ ($S_y$ resp.) after the transformation \eqref{transf}. For any element in the ring $\calO_E[X,Y,b]$ ($\calO_E[X,Y,a]$ resp.) we will denote by a bar the corresponding element in the ring $\calO_E[X,Y,b]/(bY-\pi,F_1(X,Y))$ ($\calO_E[X,Y,a]/(a\pi-Y,F_2(X,Y))$ resp.).\\ Now let $\calS,\calS^{'}$ be two horizontal divisors on $\mathfrak{F}_p$ associated with cusps $S,S^{'}$ and let $Q\in \supp \calS \cap \supp \calS^{'}$ be a point. We will denote by $m$ the maximal ideal corresponding to $Q$. If the cusps lie above different branch points, for example $S=(0:\zeta_p^{i}-1:1)$ and $S^{'}=( \zeta_p^{j}:\zeta_p^{j}-1:1)$, we have $\overline{X},\overline{X}-\overline{\zeta_p^{j}}\in m$. But then $\overline{\zeta_p^{j}}\in m$ which is impossible since $\overline{\zeta_p^{j}}$ is a unit. So  let $S$ and $S^{'}$ lie above the same branch point. Without loss of generality we may assume $S=( \zeta_p^{i}:\zeta_p^{i}-1:1)$ and $S^{'}=( \zeta_p^{j}:\zeta_p^{j}-1:1)$. It is a basic result from number theory that $(\zeta_p^i-1)/\pi$ is a unit in $\calO_E$ if $i\not\equiv0 \mod p$. We will denote this unit by $\epsilon_i$. If $Q$ is a point in the fibre $\mathfrak{F}_p \times_{\Spec \calO_E} \Spec k(\frakq)$, where $\frakq \in \Spec \calO_E$, then $\overline{\frakq} \subseteq m$. On the other hand since $\overline{X}-\overline{\zeta_p^i},\overline{X}-\overline{\zeta_p^j}\in m$ we have $\zeta_p^i-\zeta_p^j=\zeta_p^i(1-\zeta_p^{j-i})=\zeta_p^i \epsilon_{j-i} \pi$ and so $(\overline{\pi})\subseteq m$. Now if $\frakq$ is different from $(\pi)$ and so in particular coprime to $(\pi)$ we have $\overline{1}\in m$ which gives us a contradiction again. It follows that the only possibility for $Q$ to be in a special fibre is to be in the fibre of bad reduction $\mathfrak{F}_p \times_{\Spec \calO_E} \Spec k(\pi)$. Now since $S$ and $S^{'}$ are $E$-rational points $\calS$ and $\calS^{'}$ are reduced to single points $P$ and $P^{'}$ in this fibre. A direct computation shows that \[ M=\left(\overline{X}-\overline{\zeta_p^i}, \overline{\pi}, \overline{a}-\overline{\epsilon_i} \right) \] and \[M^{'}=\left(\overline{X}-\overline{\zeta_p^j}, \overline{\pi}, \overline{a}-\overline{\epsilon_j} \right)\] are the ideals corresponding to these points. If we take a look at the affine open set $U_2$, described in the previous section, we can easily verify that $M$ and $M^{'}$ are indeed maximal ideals and that $\calS$ and $\calS^{'}$ are reduced to these points in the fibre of bad reduction since \[ \overline{\pi}\left(\overline{a}-\overline{\epsilon_i}\right)=\overline{Y}-\overline{\zeta_p^i}+\overline{1} \] and $\overline{\pi}\left(\overline{a}-\overline{\epsilon_j}\right)=\overline{Y}-\overline{\zeta_p^j}+\overline{1} $. Now if $P=P^{'}=Q$ we have \[ \epsilon_i-\epsilon_j=\frac{\zeta_p^i-1}{\pi}-\frac{\zeta_p^j-1}{\pi}=\frac{\zeta_p^i-\zeta_p^j}{\pi}=\frac{\zeta_p^i(1-\zeta_p^{j-i})}{\pi}=\zeta_p^i\epsilon_{j-i} \, .\] and so $\overline{\zeta_p^i\epsilon_{j-i}}\in m$. But since $\zeta_p^i\epsilon_{j-i}\in \calO_E^{\ast}$, this gives us a contradiction and we have completed the proof of $(i)$.\\
Now let $S=(0:\zeta_p^i-1:1)$, so $S$ is $S_x$ after the transformation \eqref{transf}. Again $\calS\cap\mathfrak{F}_p \times_{\Spec\calO_E}\Spec k(\pi)$ is reduced to a single point $P$. Let $M$ be the corresponding maximal ideal, so $M=(\overline{X},\overline{\pi},\overline{a}-\overline{\epsilon_i})$. The irreducible component $L_x$ corresponds (in $U_2$) to the prime ideal $I=(\overline{\pi},\overline{X})$. Obviously $I\subset M$ and so $P$ is just in the component $L_x$ in the fibre of bad reduction (remember that the component $L$ does not lie in $U_2$). Since $\calS$ is only reduced to $P$ it only intersects $L_x$. Similar computations for $S_y$ and $S_z$ yield $(ii)$.
\end{beweis}

%
%

\begin{lemma}\label{hlemma3} Let $\mathfrak{F}_p\rightarrow \Spec\mathcal{O}_E$ be
  the arithmetic surface constructed above. There exists a canonical
  divisor $\calC\in\Cl (\mathfrak{F}_p)_\mathbb{Q} =
  \Cl(\mathfrak{F}_p)\otimes_\mathbb{Z} \mathbb{Q}$ on
  $\mathfrak{F}_p$ of the form \[ \calC=(2g-2)\calS +
    \calV \, ,\] where $\calS$ is a horizontal divisor
  coming from a cusp, $g=g(F_p)$ is the genus of $F_p$ and $\calV$
  denotes a vertical divisor having support in the special fibre
  $\mathfrak{F}_p\times_{\Spec \calO_E}{\Spec k(\pi)}$.
\end{lemma}
\begin{beweis}
It follows from Proposition \ref{prop:Rohrlich} that \[ (2g-2)S \] is a canonical divisor in $\Cl(F_p)_\QQ$, where $S$ is any cusp. If we now set \[\calC_0:=(2g-2)\calS+\calV_0 \, ,\]
  where $\calS$ is the Zariski closure of $S$ and $\calV_0$ is a sum of divisors, having support in the closed fibres, so that $\calC_0$ fulfills the adjunction formula, then $\calC_0$ is a canonical divisor of $\mathfrak{F}_p$ (see Proposition \ref{hprop1}). Note that similar arguments, as in the proof of Proposition \ref{hprop1}, assure that $\calV_0$ exists. For all primes $\frakq\in \Spec \calO_E$ not dividing $p$ - in fact these are the primes of good reduction - the special fibre $\mathfrak{F}_p\times_{\Spec \calO_E}{\Spec k(\frakq)}$ is smooth and so it consists of a single irreducible component. Since the self-intersection of this fibre is zero (see \cite{Lang}: p.61: Proposition 3.5.) we can add any multiple of it to $\calC_0$ and the resulting divisor still fulfills the adjunction formula. Using this fact we can transform $\calC_0$ into a divisor $\calC=(2g-2)\calS+\calV$, where $\calV$ is a vertical divisor having support in the special fibre over $\pi$. Again, by Proposition \ref{hprop1}, this is a canonical divisor. 
\end{beweis}

Now we are ready to compute the canonical divisor for the model
$\mathfrak{F}_p$. In the previous lemma we saw that such a divisor can
be constructed with a horizontal divisor $\calS$ coming from a cusp
and vertical divisors having support in the fibres of bad reduction. Now let $S_x$ be a cusp, \begin{equation}\label{notexyz}\calV_{x}=\lambda_x L_x
  +\lambda_y L_y +\lambda_z L_z\end{equation} and
\begin{equation}\label{notesum}\calV_{\Sigma}=\sum_{i=1}^{r} \left(
    \sum_{j=1}^{p} \lambda_{\alpha_{i,j}}L_{\alpha_{i,j}}
    +\lambda_{\alpha_i}L_{\alpha_i} \right)+\sum_{i=1}^{s
  }\lambda_{\beta_i}L_{\beta_i} \, ,\end{equation} where 
  \begin{eqnarray}
      \lambda_x &=& \left( \frac{2g-p}{p} \right) \, ,\label{l1}\\
\lambda_y=\lambda_z=\lambda_{\beta_i}=\lambda_{\alpha_{j,k}}&=&-\left(\frac{p-2}{p}\right) \mbox{ for all } i=1,\ldots,s\mbox{ and } j=1,\ldots,r \, ,\label{l2}\\
\lambda_{\alpha_j}&=&-2\left(\frac{p-2}{p} \right)\mbox{ for all } j=1,\ldots,r \, .\label{l3}                                     \end{eqnarray}
Then we claim that the divisor $\calC_x$ given by
\begin{equation} \label{abkur}\calC_x=(2g-2)\calS_x
  +\calV_{x}+\calV_{\Sigma} \, \end{equation} is a canonical divisor. Notice that $L$ is not included in $\calC_x$, since it
is  modulo the full fiber just a linear combination of the other components. 

\begin{lemma}\label{koef}
The divisor $\calC_x$ in \eqref{abkur} is indeed a canonical divisor.
\end{lemma}
\begin{beweis} From Lemma \ref{hlemma3} we know that there exists a canonical divisor of the form \eqref{abkur} with \eqref{notexyz} and \eqref{notesum} for some coefficients $\lambda$. The only thing we need to do is to show that for these $\lambda$ is no other choice possible than the one we made in \eqref{l1}, \eqref{l2} and \eqref{l3}. So the whole idea of the proof is the repeating use of the adjunction formula (see \cite{Liu}, p.390: Theorem 1.37) combined with the fact that the genus of the components of the special fibre is zero (see \cite{Mc}, p.59: Theorem 3) to approve the choice we made. We start with the observation \begin{equation}\label{alpha} 2\lambda_{\alpha_{i,j}}=\lambda_{\alpha_i}\,.\end{equation} Indeed, according to the adjunction formula $L_{\alpha_{i,j}}^2+\calC_x \cdot L_{\alpha_{i,j}}=2g(L_{\alpha_{i,j}})-2$ and $L_{\alpha_{i,j}}^2=-2$ (see previous section), we have \[ 0= \Laa \cdot\calC_x=\Laa \cdot \left( \sum_{l=1}^p \lambda_{\alpha_{i,j}}\La{i,j}+\lambda_{\alpha_i}\La{i}\right) =\lambda_{\alpha_{i,j}}(-2)+\lambda_{\alpha_i}\, .\] Now using \eqref{alpha} and the formula for $\La{i}$, we get \[ p-2= \La{i}\cdot\calC_x=\sum_{j=1}^p \lambda_{\alpha_{i,j}} +\lambda_{\alpha_{i}}(-p)=\frac{p}{2}\lambda_{\alpha_{i}}-p\lambda_{\alpha_{i}}=-\frac{p}{2}\lambda_{\alpha_{i}}\,. \] Similar computations yield $\lambda_y, \lambda_z$ and the $\lambda_{\beta_i}$. Finally, one observes that \[ p-2= \calC_x \cdot L_x=(2g-2)\calS_x  \cdot L_x +\lambda_x L_x^2 = (2g-2)+\lambda_x (-p) \] and with this we finish our proof. 
\end{beweis}

With a view to this lemma we see that the vertical part of two divisors coming from cusps that lie over different branch points, say $\calC_x$ and $\calC_y$, just differs in the parts $\calV_{x}$ and $\calV_{y}$.

\section{The algebraic contributions to \texorpdfstring{${\overline{\omega}_\Ar}^2$}{OmegaAr2}}

We now calculate certain intersection numbers, which will be used
later to complete the computations of the coefficient $a_p$.

\begin{lemma}\label{HLemma1}
  For $\calV_{\Sigma}$ given in \eqref{notesum} we have \[
  \calV_{\Sigma}\cdot \calV_{\Sigma}= (p-3)(-p)\left( \frac{p-2}{p}
  \right)^{2} \,. \]
\end{lemma}

\begin{beweis}
  In all the computations in this proof we have to remember the
  coefficients we calculated in Lemma \ref{koef}. If we write
  $\calV_{\Sigma}=\calV_{\Sigma_\alpha}+\calV_{\Sigma_\beta}$, where
  $\calV_{\Sigma_\alpha}$ denotes the part coming from the
  $L_{\alpha}$ and $\calV_{\Sigma_\beta}$ the part coming from the
  $L_{\beta}$, we have \[ \calV_{\Sigma} \cdot \calV_{\Sigma}=
  \calV_{\Sigma_\alpha}\cdot \calV_{\Sigma_\alpha} +
  \calV_{\Sigma_\beta}\cdot \calV_{\Sigma_\beta} \, ,\] since each of
  the components of $\calV_{\Sigma_\alpha}$ does not intersect any
  component of $\calV_{\Sigma_\beta}$ and vice versa. From
  figure \ref{d1} we see that each $L_{\beta_i}$ just intersects itself
  and that the number of self-intersection is $-p$. Since there are
  $s$ lines $L_{\beta_i}$, we have \[ \calV_{\Sigma_\beta}\cdot
  \calV_{\Sigma_\beta}= s  (-p) \left(\frac{p-2}{p}
  \right)^{2}  \, .\] Now let $\calC$ be a canonical
  divisor. According to the adjunction formula, we have $\calC\cdot
  L_{\alpha_{i,j}}=0$ and, since each $L_{\alpha_{i,j}}$ just
  intersects the $\calV_{\Sigma_{\alpha}}$ part of $\calC$, the
  equation $0=\calC\cdot L_{\alpha_{i,j}}= \calV_{\Sigma_\alpha} \cdot
  L_{\alpha_{i,j}}$. This 
yields \[ \calV_{\Sigma_\alpha} \cdot \calV_{\Sigma_\alpha}
=\calV_{\Sigma_\alpha} \cdot \sum_{i=1}^{r} \lambda_{\alpha_{i}}
L_{\alpha_i} = \sum_{i=1}^{r}
\lambda_{\alpha_{i}}\left(\calV_{\Sigma_\alpha} \cdot
  L_{\alpha_i}\right) \, ,\] where each addend is 
\begin{eqnarray*} \lambda_{\alpha_i}\left(\calV_{\Sigma_\alpha} \cdot L_{\alpha_i}\right)&=& \lambda_{\alpha_i} \left(\left( \sum_{i=1}^p \lambda_{\alpha_{i,j}} L_{\alpha_{i,j}} + \lambda_{\alpha_i} L_{\alpha_i} \right) \cdot L_{\alpha_i} \right)  \\ &=&\lambda_{\alpha_i} \left( \frac{p}{2} \lambda_{\alpha_i} + \lambda_{\alpha_i} (-p) \right) \\
  &=& -\frac{p}{2}\lambda_{\alpha_i}^2 =
  2(-p)\left(\frac{p-2}{p}\right)^2 \, .\end{eqnarray*}
Since there are $r$ lines $L_{\alpha_i}$, we have \[
\calV_{\Sigma}\cdot \calV_{\Sigma}= (2r+s) (-p)\left( \frac{p-2}{p}
\right)^{2}= (p-3)(-p)\left( \frac{p-2}{p} \right)^{2} \]
\end{beweis}
\begin{lemma}\label{HLemma2}
 Let $\calV_x$ be a vertical divisors as in \eqref{notexyz} which belongs to a cusp. Then \[
                                                         \calV_{x}\cdot \calV_{x}=(-p)\left(\frac{2g-p}{p}\right)^2  +(-2p)\left(\frac{p-2}{p}\right)^{2} .\]
\end{lemma}
\begin{beweis}
The lines $L_x$,$L_y$ and $L_z$ only intersect themselves and each self-intersection number is $-p$. Now everything follows from the equations \eqref{l1} and \eqref{l2}.
\end{beweis}

\begin{lemma}
Let \begin{equation}\label{o1} \calD_x =\calS_x +\calG_x \, ,
\end{equation}where $\calG_x=\frac{1}{p}L_x$. Then the divisor $\calD_x$ is associated with $\left(\beta^{\ast}\calO_{\mathbb{P}_{\calO_E}^1}(1)\right)^{\otimes \frac{1}{p^2}}$, or in other words $\calO (D_x)^{\otimes p^2}\cong \beta^{\ast}\calO_{\mathbb{P}_{\calO_E}^1}(1)$.
\end{lemma}

\begin{beweis}
Let $S_x$ be a cusp and $Q\in\PP_E^1$ the corresponding branch point. Since $\Pic(\PP_E^1)\cong \ZZ$ and $\calO_{\PP_E^1}(1)$ is a generator of $\Pic(\PP_E^1)$ any divisor of degree $1$ is associated with $\calO_{\PP_E^1} (1)$. We choose $Q$ to be this associated divisor. Now \[ \beta^{\ast}Q=\sum_{i=1}^{p} p S_i \, ,\] where $S_i$ runs through the cusps lying above $Q$. If follows from \cite{Roh}, p.101: Theorem 1. that $\beta^{\ast} Q \sim p^2  S_x$ in $\Cl(F_p)_\QQ$ (remember that $S_x$ is one of the cusps) and so $p^2  S_x$ is associated with $\beta^{\ast}\calO_{\PP_E^1} (1)$. Since $\beta^{\ast} \calO_{\PP_{\calO_E}^1} (1)|_{F_p} \sim \beta^{\ast}\calO_{\PP_E^1} (1)$ it is clear with Lemma \ref{hlemma2} that we can choose $\calD_x=\calS_x + \calG_x$ where $\calG_x$ is a vertical divisor having support in the special fibre $\mathfrak{F}_p \times_{\Spec \calO_E} \Spec k(\pi)$. Now let $I$ be the index set from \eqref{indexset}. Since each component of the special fibre which is different to $L$ is mapped to a single point by $\beta$, we have \begin{equation}\label{e1}
(p^2 \calD_x)\cdot L_i =0 ~~~~~~~(\forall i\in I) 
\end{equation}(see \cite{Liu}, p. 398: Theorem 2.12 (a) ). On the other hand we have \begin{equation}\label{e2} p^2=p^2\calD_x \cdot \mathfrak{F}_p \times_{\Spec \calO_E} \Spec k(\pi) = p^2\calD_x \cdot L 
\end{equation} (see \cite{Liu}, p. 388: Remark 1.31.). Solving \eqref{e1} and \eqref{e2} we get $\calG_x =\frac{1}{p}L_x$.
\end{beweis}

\begin{satz}\label{thm1}
Let $\calC_x=(2g-2)(\calS_x+\calF_x)$ be a canonical divisors and $\calD_x=\calS_x+\calG_x$ a divisors as in \eqref{o1}, where $x$ indicates that this divisor belongs to a cusp $S_x$. Then \[\renewcommand{\arraystretch}{1.75}\begin{array}{ccc}
 \calF_x \cdot \calF_x  = -\frac{p^3-7p^2 +15p -8}{p^2 (p-3)^2 }\, ,     & & \calS_x \cdot \calG_x  = -(\calG_x \cdot \calG_x)                 = \frac{1}{p}  \, .                 
\end{array}\]
\end{satz}
\begin{beweis}
We have $\calF_x^2 =\frac{1}{(2g-2)^2}\left(\calV_{x}^2+\calV_{\Sigma}^2\right)$. Now Lemma \ref{HLemma1} and Lemma \ref{HLemma2} together with $g=\frac{(p-1)(p-2)}{2}$ yield (after simplifying equations) our first claim. \\ With equation \eqref{e1} we get $\calS_x \cdot \calG_x  = -(\calG_x \cdot \calG_x)$. Since $\calG_x=\frac{1}{p}L_x$ the second claim follows.
\end{beweis}


Now, we successfully prepared all the ingredients to actually calculate some intersection numbers for the Fermat curves.

\section{Proof of the main result}

\begin{satz}\label{thm:fermat}
  Let $\mathfrak{F}_{p}$ be the regular model of the fermat curve
  $F_p$ over $\Spec \ZZ[\zeta_p]$ which was constructed in section
  \ref{regandminmodel}.  Then the arithmetic self-intersection number
  of its dualizing sheaf equipped with the Arakelov metric satisfies
\begin{equation*} \overline{\omega}_{\mathfrak{F}_p, \Ar }^2 \le (2 g-2) \left( \log|\Delta_{\QQ(\zeta_p)|\QQ}|^2 + 
   [\QQ(\zeta_p):\QQ]  \left(  \kappa_1 \log p + \kappa_2 \right) + \frac{p^2-4p+2}{p(p-3)} \log p
   \right),
\end{equation*}
where $\kappa_1, \kappa_2\in \RR$ are  constants independent of $p$.
\end{satz}

\begin{beweis}
  In section \ref{sec:belyimorph} and \ref{regandminmodel} we saw that
  the morphism $\beta:F_p \rightarrow \PP^1$ fulfills the requirements
  of Theorem \ref{thm:keyformula}. Since $\beta^{\ast}\infty
  =\sum_{i=1}^p pS_i$ we have $b_j=b_{\max}=p$. The morphism $\beta$
  is of degree $p^2$. It follows that in our case the formula
  \eqref{eq:def-ap} of Theorem \ref{thm:keyformula}
  becomes \begin{align*} \sum_{\frakp \,\,{\rm bad}} a_\frakp
    \log\Nm(\frakp) = a_\frakp \log\Nm(\frakp)
    &=- 2 g\calO( \calG_j)^2 + (2g-2)\calO(\calF_j)^2 \\
    &= -2 g {\calG_j}^2 \log p + (2g-2) {\calF_j}^2 \log p \\[6pt]
    &= \frac{2 g}{p}\log  p - (2g-2) \frac{p^3-7p^2 +15p -8}{p^2 (p-3)^2 }\log p \\
    &= \frac{p^2-4p+2}{p(p-3)} \log p .
\end{align*}
\end{beweis}


\begin{bem}
In Section \ref{regandminmodel} we have seen that we get a minimal regular model $\mathfrak{F}_p^{min}$ of $F_p$ if we blow down the component $L$ of the special fibre. Let $\pi:\mathfrak{F}_p \rightarrow \mathfrak{F}_p^{min}$ denote this blow-down. Then there exists a vertical divisor $\calW$ on $\mathfrak{F}_p$ (with support in the special fibre) such that $\pi^{\ast}\omega_{\mathfrak{F}_p^{min}}=\omega_{\mathfrak{F}_p}\otimes\calO (\calW)$. We have \begin{equation*}\label{umrechnung}
\overline{\omega}_{\mathfrak{F}_p^{min} , Ar}^2= \pi^{\ast} \overline{\omega}_{\mathfrak{F}_p^{min} , Ar}^2=\overline{\omega}_{\mathfrak{F}_p, Ar}^2 +2\omega_{\mathfrak{F}_p}\cdot\calO (\calW) +\calO (\calW)^2 \, .
\end{equation*}
\end{bem}
        
\begin{proposition}\label{prop:change}
With the notation from above we have \begin{equation*}\label{umrechnung2}
 2\omega_{\mathfrak{F}_p}\cdot\calO (\calW) +\calO (\calW)^2 =(2p^2-10p+13)\log p .
\end{equation*}
\end{proposition}

\begin{beweis}
We start by computing the canonical divisor $\calK_x^{min}$ of $\mathfrak{F}_p^{min}$, so the divisor with $\calO (\calK_x^{min}) \cong \omega_{\mathfrak{F}_p^{min}}$. Let $\tilde{L}_u:=\pi L_u$, where $u\in I$ and $I$ is the index set \eqref{indexset}. In order to compute intersections of the $\tilde{L}_u$ we need to find their pullback and then compute everything on $\mathfrak{F}_p$. We have $\pi^{\ast} \tilde{L}_u =L_u$ for $u=\alpha_{i,j}$ and \[ \pi^{\ast}\tilde{L}_u=L_u + L \] for all other $u$. Indeed, let for instance $u=x$. Then we have $\pi^{\ast}\tilde{L}_x=L_x +\mu_x L$, where $\mu_x$ is a rational number. It follows that $0=L\cdot \pi^{\ast} \tilde{L}_x=1-\mu_x$ (see \cite{Liu}, p.398: Theorem 2.12. (a)). \\
The canonical divisor on $\mathfrak{F}_p^{min}$ is given by \begin{equation*}
\calK_x^{min}=(2g-2)(\calS_x +\frac{1}{p}\tilde{L}_x) \, .
\end{equation*}To verify this we just need to proof that $\calK_x^{min}$ satisfies the adjunction formula and restricts to the canonical divisor $K_x$ of the generic fibre $F_p$ (see Proposition \ref{hprop1}). The second property is obviously fulfilled. In order to verify the adjunction formula one has to check that it is valid for each irreducible component of the special fibre. We will illustrate this for the component $\tilde{L}_x$ and leave the rest to the reader since the computations are very similar. We have \begin{eqnarray*} \calK_x^{min}\cdot \tilde{L}_x &=& (2g-2)(\calS_x\cdot \tilde{L}_x +\frac{1}{p}\tilde{L}_x^2) \\ &=& (2g-2)(1+\frac{1}{p}(L_x+L)^2) \\ &=& p(p-3)(1-\frac{1}{p}(p-1)) = (p-3)\end{eqnarray*} (see \cite{Liu}, p.398: Theorem 2.12. (c) for the second equality). On the other hand is \[2p_a (\tilde{L}_x)-2-\tilde{L}_x^2=-2-(L_x+L)^2=(p-3) \] and so the formula is valid for $\tilde{L}_x$.\\
The pullback of the canonical divisor is now \[\pi^{\ast}\calK_x^{min}=(2g-2)(\calS_x + \frac{1}{p}L_x +\frac{1}{p} L) \] and an easy computation shows that \[ \calW=-\lambda_y L_y -\lambda_z L_z -\frac{(2-p)}{p} L_x -\calV_{\Sigma} +\frac{2g-2}{p} L \] fulfills $\pi^{\ast}\calK_x^{min}=\calK_x +\calW$. It follows that we have to compute $(2\calK_x \cdot \calW + \calW^2)\log p$ in order to get $2\omega_{\mathfrak{F}_p}\cdot\calO (\calW) +\calO (\calW)^2$. Since we have $\calW\cdot(2\calK_x +\calW)=\calW\cdot (\calK_x+\pi^{\ast}\calK_x^{min})$ we may compute $\calW\cdot\calK_x$ and $\calW \cdot \pi^{\ast}\calK_x^{min}$. Using the adjunction formula and linearity we get \begin{eqnarray*}
 \calW \cdot \calK_x &=& (p-2)\left(-\lambda_y-\lambda_z-\left(\frac{2-p}{p}\right) \right)-\calV_{\Sigma}\cdot \calK_x -\left(\frac{2g-2}{p} \right) \\ &=& 3 \left( \frac{(p-2)^2}{p}\right) - \calV_{\Sigma}^2 -\left(\frac{p(p-3)}{p}\right) \\ &=& (p-2)^2-(p-3) \, .
\end{eqnarray*}
On the other hand we have \begin{eqnarray*}
\calW\cdot \pi^{\ast}\calK_x^{min} &=& \calW\cdot(p(p-3)\calS_x +(p-3)L_x +(p-3)L) \\
&=& (p-2)(p-3)-(p-2)(p-3)+(p-3)^2+ (p-3)\calW \cdot L \\ 
&=& (p-3)^2+(p-3)\left(-\lambda_y-\lambda_z-\frac{2-p}{p}+\frac{p-2}{p}(p-3) -(p-3) \right) \\ 
&=& (p-3)^2+(p-3)(p-2)-(p-3)^2=(p-2)(p-3) 
\end{eqnarray*} and so $2\omega_{\mathfrak{F}_p}\cdot\calO (\calW) +\calO (\calW)^2= (2p^2-10p+13)\log p$.
\end{beweis}

\begin{satz}\label{thm:omega-min}
  Let $\mathfrak{F}_{p}^{min}$ be the minimal regular model of the
  fermat curve $F_p$ over $\Spec \ZZ[\zeta_p]$ from section
  \ref{regandminmodel}.  Then the arithmetic self-intersection number
  of its dualizing sheaf equipped with the Arakelov metric satisfies
\begin{equation*} \overline{\omega}_{\mathfrak{F}_p^{min}, \Ar }^2 \le (2 g-2) \left( \log|\Delta_{\QQ(\zeta_p)|\QQ}|^2 + 
   [\QQ(\zeta_p):\QQ]  \left(  \kappa_1 \log p + \kappa_2 \right) + \frac{3p^2-14p+15}{p(p-3)} \log p
   \right),
\end{equation*}
where $\kappa_1, \kappa_2\in \RR$ are  constants independent of $p$.
\end{satz}

\begin{beweis}
Follows directly from Theorem \ref{thm:fermat} and Proposition \ref{prop:change}.
\end{beweis}

\begin{korollar}
With the notation from the previous theorem we have:
\begin{equation*} \overline{\omega}_{\mathfrak{F}_p^{min}, \Ar }^2 \le (2 g-2) \left(  
   (p-1)  \left(  \kappa_1 \log p + \kappa_2 \right) + \frac{2p^2-p-5}{p} \log p \right)
\end{equation*}
\end{korollar}
\begin{beweis}
  It is a well known fact that
  $\Delta_{\QQ(\zeta_p)|\QQ}=(-1)^{\frac{p-1}{2}}p^{p-2}$ and
  $[\QQ(\zeta_p):\QQ]=p-1$ and so Theorem \ref{thm:omega-min}
  yields \begin{eqnarray*} \overline{\omega}_{\mathfrak{F}_p, \Ar
    }^2&\le & (2 g-2) \left( \log p^{2p-4} +
      (p-1)  \left(  \kappa_1 \log p + \kappa_2 \right) + \frac{3p^2-14p+15}{p(p-3)} \log p \right) \\
    &=& (2 g-2) \left( (p-1) \left( \kappa_1 \log p + \kappa_2 \right)
      + \frac{2p^2-p-5}{p} \log p\right)
\end{eqnarray*}
\end{beweis}


%

\def\cprime{$'$}
\providecommand{\arxivref}[1]{\href{http://arxiv.org/abs/#1}{#1}}
\providecommand{\bysame}{\leavevmode\hbox to3em{\hrulefill}\thinspace}
\providecommand{\MR}{\relax\ifhmode\unskip\space\fi MR }
\providecommand{\MRhref}[2]{%
  \href{http://www.ams.org/mathscinet-getitem?mr=#1}{#2}
}
\providecommand{\href}[2]{#2}

\end{document}